\documentclass[11pt]{article}
\usepackage{amssymb,latexsym,color,amsmath,epsfig}
\usepackage[top=0.8in,bottom=1in,left=1in,right=1in]{geometry}
 \usepackage{graphicx}   
\usepackage{epsfig}
\usepackage{graphicx,graphics}
\usepackage{mathrsfs}
\usepackage{times}
\usepackage{mathptmx}
\usepackage{amssymb,amsmath,amsthm}
\usepackage{amsfonts}
\usepackage{txfonts}
\usepackage{subfigure}
\usepackage{amsfonts} \usepackage{txfonts}
\newtheorem{theorem}{\textbf {Theorem}}
\newtheorem{assumption}{\textbf {Assumption}}
\newtheorem{proposition}{\textbf {Proposition}}
\newtheorem{lemma}{\textbf {Lemma}}

\newtheorem{defn}{\textbf {Definition}}

\begin{document}
\title{A Broad Dynamical Model for Pattern Formation by Lateral Inhibition}
\author{Murat Arcak\thanks{Department of Electrical Engineering and Computer Sciences,
University of California, Berkeley. Email:  arcak@eecs.berkeley.edu.}}
\date{\today}
\maketitle

\section{Introduction}
Spatial patterns of gene expression are central to the development of multi-cellular organisms. Most mathematical studies of pattern formation investigate {\it diffusion-driven instability}, which is a mechanism that amplifies spatial inhomogeneities in a class of reaction-diffusion systems (see, {\it e.g.}, \cite{murray2}). However, many patterning events in multi-cellular organisms rely on cell-to-cell contact signaling, such as the {\it Notch pathway} \cite{gilbert}, and do not involve diffusible proteins for intercellular communication.  A particularly interesting phenomenon in this form of communication is
{\it lateral inhibition} whereby a cell that adopts a particular fate inhibits its immediate neighbors from doing likewise \cite{ColMonMaiLew96}, thus leading to `fine-grained' patterns. There is increasing interest in understanding the Notch signaling circuitry in mammalian cells that leads to such lateral inhibition \cite{SprLakLeB10,SprLakLeB11}. Recent studies showed that a  lateral inhibition pathway also functions in {\it E. Coli}, and enables the bacteria to inhibit the growth of other {\it E. Coli} strains in direct contact \cite{AokLow10}.

Dynamical models are of great interest for understanding the circuit topologies involved in lateral inhibition and for predicting the associated patterns. Several simplified models have been employed for Notch signalling pathways in \cite{ColMonMaiLew96} and \cite{SprLakLeB11}. The objective of this paper is to present an abstract dynamical model that captures the essential features of lateral inhibition and to demonstrate with dynamical systems techniques that these features indeed lead to patterning. Although this model is not meant specifically for Notch signaling, it encompasses as special cases the lateral inhibition model in \cite{ColMonMaiLew96} as well as a slightly modified version of the one in \cite{SprLakLeB11}.

Our model treats the evolution of concentrations in each cell as an input-output system, where the inputs represent the influence of adjacent cells and the outputs correspond to the concentrations of the species that interact with adjacent cells. The input-output models for the cells are then interconnected according to an undirected graph where the nodes represent the cells, and the presence of a link between two nodes means that the corresponding cells are in contact.  The main assumption on the input-output model is that each constant input yields a unique and globally asymptotically stable steady-state, and that the value of the output at this steady-state is a decreasing function of the input. This decreasing property captures the inhibition of the cell function by its neighbors.
The model allows multiple inputs and outputs, and is restricted by a {\it monotonicity} assumption, following the definition of monotonicity for dynamical systems with inputs and outputs \cite{AngSon03}.

Using this model, we first give an instability condition for the homogeneous steady-state, applicable to arbitrary contact graphs.  We then focus our attention on bipartite graphs, and demonstrate the emergence of a ``checkerboard" pattern, exhibiting alternating high and low values of concentrations in adjacent cells. Next, we establish a {\it strong monotonicity} property  of the interconnected model for bipartite graphs, which implies that almost every bounded solution (except for a measure-zero set of initial conditions) converges to a steady-state \cite{Hirsch-Smith,smith}. A graph is bipartite if and only if it contains
 no odd-length cycles, and Cartesian products of bipartite graphs are also bipartite \cite{asratian}.
Thus, the results of this section are applicable, among others, to {\it grid graphs} (one dimensional {\it path graphs} and their Cartesian products in higher dimensions) which are appropriate for representing arrays of cells.

\section{Lateral Inhibition Model and Preliminaries}\label{CIM}
We let $\mathcal{G}$ be an undirected, connected graph where the nodes represent the cells, and the presence of a link between two nodes means that the corresponding cells are in contact. In preparation for the dynamical model studied below, we let $N$ denote the number of cells and define the matrix
 $P\in \mathbb{R}^{N\times N}$:
\begin{equation}\label{randomwalk}
p_{ij}=\left\{ \begin{array}{ll}{d_i}^{-1} & \mbox{if nodes $i$ and $j$ are adjacent,} \\ 0 & \mbox{otherwise,}\end{array}\right.
\end{equation}
where $d_i$ denotes the degree of node $i$.  It follows that $P$ is a nonnegative {\it row-stochastic} matrix, that is:
\begin{equation}
P\bf{1}=\bf{1} \label{ones}
\end{equation}
where $\bf{1}$ denotes the vector of ones.
The matrix $P$ is identical to the probability transition matrix for a {\it random walk} on the graph $\mathcal{G}$.  The properties summarized below therefore follow from standard results for random walks (see, {\it e.g.}, \cite{levin}):

\begin{lemma} \label{prop}
$P$ possesses real eigenvalues $\lambda_N\le \cdots \le \lambda_1$ all of which lie in the interval $[-1,1]$, and corresponding real, linearly independent eigenvectors $v_i$, $i=1,\cdots,N$.
In particular, $\lambda_1=1$, and $v_1=\mathbf{1}$ is a corresponding eigenvector. If $\mathcal{G}$ is bipartite, then $\lambda_N=-1$, and an eigenvector $v_N$ is such that the entries are either $1$ or $-1$, and two entries corresponding to adjacent nodes have opposite signs.
\end{lemma}

Let $i=1,\cdots,N$ denote the cells, and consider the dynamical model:
\begin{equation}
\label{one}
\dot{x}^i=f(x^i,u^i) \quad y^i=h(x^i)
\end{equation}
where $x^i\in \mathscr{X}\subset\mathbb{R}^n$ is a vector describing the state of reagent concentrations in cell $i$, $u^i\in \mathscr{U}\subset \mathbb{R}^m$ describes the `input' from adjacent cells, and $y^i\in \mathscr{Y}\subset\mathbb{R}^m$ describes the `output' that serves as an input to adjacent cells. In particular,
\begin{equation}
\label{two}
U=(P\otimes I_m)Y
\end{equation}
where $P$ is as defined in (\ref{randomwalk}),
$U:=[{u^1}^T \cdots {u^N}^T]^T$ and $Y:=[{y^1}^T \cdots {y^N}^T]^T$. If follows from (\ref{randomwalk}) that the input $u^i$  is the average of the outputs $y^k$ over all neighbors $k$ of cell $i$. Thus, we henceforth take the input and output spaces to be identical: $\mathscr{U}=\mathscr{Y}$.

We assume that $f(\cdot,\cdot)$ and $h(\cdot)$ are continuously differentiable and further satisfy the following property:
\begin{assumption}
\label{characteristic}
For each constant input $u^*$, system (\ref{one}) has a globally asymptotically stable steady-state $x^*:=S(u^*)$ with the additional property that:
\begin{equation}\label{det}
\det\left(\left.\frac{\partial f(x,u)}{\partial x}\right|_{(x,u)=(x^*,u^*)}\right)\neq 0.
\end{equation}
The map $S:\mathscr{U} \rightarrow \mathscr{X}$ and, therefore, the map $T:\mathscr{U} \rightarrow \mathscr{U}$ defined by:
\begin{equation}\label{Tdef}
T(\cdot):=h(S(\cdot)),
\end{equation}
are continuously differentiable.
\end{assumption}

Following the terminology in \cite{AngSon03}, we will refer to $S(\cdot)$ as the {\it input-state characteristic}, and to $T(\cdot)$ as the {\it input-output characteristic.} Our next assumption is that (\ref{one}) is a {\it monotone} system in the sense of \cite{AngSon03}, as defined below. According to the classical definition for systems without inputs and outputs  \cite{smith}, a monotone system is one
 that preserves a partial ordering of the initial conditions as the solutions evolve.
The partial ordering is defined with respect to a {\it positivity cone} $K$ in the Eucledean space that is  closed, convex, {\it pointed} ($K\cap (-K)=\{0\}$), and has nonempty interior. Given such a cone, $x\preceq \hat{x}$ means $\hat{x}-x\in K$, $x\prec \hat{x}$ means $x\preceq \hat{x}$ and $x\neq \hat{x}$, and $x\ll \hat{x}$ means
that $\hat{x}-x$ is in the interior of $K$.
The system $\dot{x}=f(x)$ is then defined to be {\it monotone} if two solutions $x(t)$ and $\hat{x}(t)$ starting with the order $x(0)\preceq \hat{x}(0)$ maintain $x(t)\preceq \hat{x}(t)$ for all\footnote{Here, ``for all $t$" is understood as ``for all times $t$ that belong to the common domain of existence  of the two solutions."} $t\ge 0$. The more restrictive notion of {\it strong monotonicity} stipulates that $x(0)\prec \hat{x}(0)$  implies $x(t)\ll \hat{x}(t)$ for all $t>0$. The monotonicity concept was extended to systems with inputs and outputs in \cite{AngSon03}:

\begin{defn}\label{angson}
Given positivity cones $K^U,K^Y,K^X$ for the input, output, and state spaces,  the system $\dot{x}=f(x,u)$, $y=h(x)$ is said to be monotone if $x(0)\preceq \hat{x}(0)$ and $u(t)\preceq \hat{u}(t)$  for all $t\ge 0$ imply that the resulting solutions satisfy $x(t)\preceq \hat{x}(t)$ for all $t\ge 0$, and the output map is such that $x\preceq \hat{x}$ implies $h(x)\preceq h(\hat{x})$.
\end{defn}

\begin{assumption}\label{iomon}
The system (\ref{one}) is monotone with respect to $K^U=\mathbb{R}_{\ge 0}^m$, $K^Y=-K^U$, and $K^X=K$, where $K$ is some positivity cone in $\mathbb{R}^n$.
\end{assumption}

As observed in \cite[Remark V.2]{AngSon03}, monotonicity implies that the input-state and input-output characteristics are nondecreasing with respect to the same ordering; that is, $u \preceq \hat{u}$ with respect to $K^U$ implies $S(u) \preceq S(\hat{u})$ with respect to $K^X$ and $T(u) \preceq T(\hat{u})$ with respect to $K^Y$. Since $K^Y=-K^U$ in Assumption \ref{iomon}, we conclude that $T(\cdot)$ is {\it nonincreasing} with respect to the standard order induced by $K^U=\mathbb{R}_{\ge 0}^m$.  This nonincreasing property means that, if two cells are in contact, an increase in the output value of one has the opposite effect on the other, which is why (\ref{one})-(\ref{two}) is referred to as a  ``lateral inhibition" model.
We note from the nonincreasing property of $T(\cdot)$ that:
\begin{equation}\label{Tprime}
T'(u):=\frac{\partial T(u)}{\partial u}
\end{equation}
is a nonpositive matrix in $\mathbb{R}^{m\times m}$, and denote its {\it spectral radius} as:
\begin{equation}
\rho(T'(u)).
\end{equation}
We conclude this section by quoting lemmas that will be used in the sequel. Lemmas \ref{jaclin} and \ref{linchar} are from \cite{AngSon04}:
\begin{lemma}\label{jaclin}
Given the system $\dot{x}=f(x,u)$, $y=h(x)$ with continuously differentiable $f(\cdot,\cdot)$ and $h(\cdot)$, the linearization $\dot{x}=Ax+Bu$, $y=Cx$ about a point $(x^*,u^*)$ satisfying $f(x^*,u^*)=0$ is also monotone with respect to the same positivity cones.
\end{lemma}
\begin{lemma} \label{linchar}
The linear system $\dot{x}=Ax+Bu$, $y=Cx$ is monotone if and only if:

1) $x\in K^X$ implies $Ax\in K^X$,

2) $u\in K^U$ implies $Bu\in K^X$,

3) $x\in K^X$ implies $Cx\in K^Y$.
\end{lemma}
The following lemma, proven in \cite{AngSon04} for single-input, single-output systems and extended in \cite{EncSon05} to the multivariable case, determines stability of a positive feedback system based on the `dc gain' of the open-loop system:
\begin{lemma} \label{posfbk}
Suppose the linear system $\dot{x}=Ax+Bu$, $y=Cx$ is monotone with respect to cones $K^U,K^Y,K^X$ such that $K^U=K^Y$ and $A$ is Hurwitz. If $-(I+CA^{-1}B)$ is Hurwitz, then so is $A+BC$. If $-(I+CA^{-1}B)$ has an eigenvalue with a positive real part, then so does $A+BC$.
\end{lemma}
\noindent
In the special case of single-input, single-output systems, the stability condition above amounts to checking whether the dc gain $-CA^{-1}B$ is greater or smaller than one.  In the multi-input, multi-output case, this condition is equivalent to inspecting whether the spectral radius of the dc gain matrix is greater or smaller than one.

The following test from \cite{AngSon03,AngSon04b} is useful for certifying monotonicity with respect to orthant cones:
\begin{lemma}\label{monotonetest} Consider the system $\dot{x}=f(x,u)$, $y=h(x)$, $x\in \mathscr{X}\subset \mathbb{R}^n$,  $u\in \mathscr{U}\subset \mathbb{R}^m$, $y\in \mathscr{Y}\subset \mathbb{R}^p$, where the interiors of $\mathscr{X}$ and $\mathscr{U}$ are convex, and
$f(\cdot,\cdot)$ and $h(\cdot)$ are continuously differentiable. If there exist $\epsilon_1,\cdots,\epsilon_n,\delta_1,\cdots,\delta_m,\mu_1,\cdots,\mu_p \in \{0,1\}$ such that:
\begin{eqnarray}
&&(-1)^{\epsilon_j+\epsilon_k}\frac{\partial f_j}{\partial x_k}(x,u)\ge 0 \quad \forall x\in \mathscr{X}, \forall u\in \mathscr{U}, \forall j\neq k\\
&&(-1)^{\epsilon_j+\delta_k}\frac{\partial f_j}{\partial u_k}(x,u)\ge 0 \quad \forall x\in \mathscr{X}, \forall u\in \mathscr{U}, \forall j,k\\
&&(-1)^{\epsilon_j+\mu_k}\frac{\partial h_k}{\partial x_j}(x,u)\ge 0 \quad \forall x\in \mathscr{X},  \forall j,k,
\end{eqnarray}
then the system is monotone with respect to the positivity cones $K^U=\{u\in \mathbb{R}^m \ | \ (-1)^{\delta_j}u_j\ge 0\}$, $K^X=\{x\in \mathbb{R}^n \ | \ (-1)^{\epsilon_j}x_j\ge 0\}$, $K^Y=\{y\in \mathbb{R}^p \ | \ (-1)^{\mu_j}y_j\ge 0\}$.
\end{lemma}

\section{Instability of the Homogeneous Steady-State}\label{INS}
Note that system (\ref{one})-(\ref{two}) admits spatially homogeneous solutions of the form
 $x^i(t)={\mathbf x}(t)$, $i=1,\cdots,N$, where ${\mathbf x}(t)$ satisfies:
\begin{equation}
\dot{{\mathbf x}}= f({\mathbf x},h({\mathbf x})).\label{single}
\end{equation}
In  particular, if the map $T(\cdot)$ has a fixed point:
\begin{equation}\label{refo1}
{\mathbf u}^*=T({\mathbf u}^*),
\end{equation}
then (\ref{single}) admits the steady-state:
\begin{equation}\label{refo2}
{\mathbf x}^*=S({\mathbf u}^*).
\end{equation}
For single-input, single-output systems with $\mathscr{U}=\mathbb{R}_{\ge 0}$, the nonincreasing property of the map $T:\mathbb{R}_{\ge 0}\rightarrow \mathbb{R}_{\ge 0}$ indeed guarantees a unique fixed point ${\mathbf u}^*$ in (\ref{refo1}).

The ``lumped model" (\ref{single}) describes the dynamics of the $Nn$-dimensional system (\ref{one}) reduced to the $n$-dimensional invariant subspace where the solutions are spatially homogeneous. Thus, the steady-state ${\mathbf x}^*$ of the lumped model defines the homogeneous steady-state $x^i={\mathbf x}^*$, $i=1,\cdots,N$, for the full system (\ref{one})-(\ref{two}). As a starting point for the analysis of pattern formation, we now give an instability criterion for the homogeneous steady-state:

\begin{theorem} \label{T1} Consider the system (\ref{one})-(\ref{two}) and suppose Assumptions \ref{characteristic} and \ref{iomon} hold. Let $\lambda_N$ denote the smallest eigenvalue of $P$ as in Lemma \ref{prop}, and let ${\mathbf u}^*$, ${\mathbf x}^*$ be as in (\ref{refo1}), (\ref{refo2}).
If:
\begin{equation}\label{unstab}
\lambda_N \, \rho\left(T'({\mathbf u}^*)\right)<-1,
\end{equation}
then the homogeneous steady-state $x^i={\mathbf x}^*$, $i=1,\cdots,N,$
is unstable.
 \end{theorem}

\noindent
{\it Proof:} Let $X:=[{x^1}^T \cdots {x^N}^T]^T$, and note that
the linearization of (\ref{one})-(\ref{two}) about the homogeneous steady-state $[{{\mathbf x}^*}^T, \cdots,{{\mathbf x}^*}^T]^T$ gives the Jacobian matrix:
\begin{equation}\label{jaco}
I_N\otimes A+P\otimes (BC)
\end{equation}
where:
\begin{equation}\label{abcdef}  A:=\left.\frac{\partial f(x,u)}{\partial x}\right|_{(x,u)=({\mathbf x}^*,{\mathbf u}^*)}, \quad B:=\left.\frac{\partial f(x,u)}{\partial u}\right|_{(x,u)=({\mathbf x}^*,{\mathbf u}^*)}, \quad C:=\left.\frac{\partial h(x)}{\partial x}\right|_{x={\mathbf x}^*}.
\end{equation}
We recall  from Lemma \ref{prop} that
\begin{equation}
V^{-1}PV=\Lambda:=\left[\begin{array}{ccc} \lambda_1 & & \\ & \ddots & \\ & & \lambda_N \end{array}\right],
\end{equation}
where $V=[v_1\cdots v_N]$, and apply the following similarity transformation to (\ref{jaco}):
\begin{equation}\label{similarity}
(V^{-1}\otimes I_n)[I_N\otimes A+P\otimes (BC)](V\otimes I_n)=I_N\otimes A+\Lambda \otimes (BC).
\end{equation}
This matrix is block-diagonal, with the $k$th diagonal block given by:
\begin{equation}\label{eachk}
A+\lambda_k BC.
\end{equation}
\noindent
{\it Claim: If
\begin{equation}\label{unstabk}
\lambda_k \, \rho\left(T'({\mathbf u}^*)\right)<-1,
\end{equation}
then (\ref{eachk}) has a positive eigenvalue.}

The theorem follows from this claim because, if (\ref{unstab}) holds, then (\ref{eachk}) has a positive eigenvalue for $k=N$, which implies instability. To prove the claim, we note from Assumption \ref{iomon} and Lemma \ref{jaclin} that the linear system $\dot{x}=Ax+Bu$, $y=Cx$ is monotone with respect to $K^U=\mathbb{R}_{\ge 0}^m$, $K^Y=-K^U$, and $K^X=K$. We write $A+\lambda_kBC=A+B{C}_k$ where ${C}_k:=\lambda_kC$ and note that (\ref{unstabk}) implies $\lambda_k<0$. Thus, the  system $\dot{x}=Ax+Bu$, $y=C_kx$ is monotone with $K^U=K^Y$. In addition, Assumptions \ref{characteristic} and \ref{iomon} imply that $A$ is Hurwitz, as can be deduced from \cite[Lemma 6.5]{AngSon04}. Thus, it follows from the second statement of Lemma \ref{posfbk} that if $-(I+C_kA^{-1}B)$ has a positive eigenvalue, then so does (\ref{eachk}).
The remaining task is thus to prove that
\begin{equation}\label{stepby}
-(I+C_kA^{-1}B)=-I-\lambda_kCA^{-1}B
\end{equation}
has a positive eigenvalue. To this end, we first show that
\begin{equation}\label{acdc}
T'({\mathbf u}^*)=-CA^{-1}B.
\end{equation}
Since
\begin{equation}
f(S(u),u)\equiv 0,
\end{equation}
differentiation gives:
\begin{equation}\label{x}
\left.\frac{\partial f(x,u)}{\partial x}\right|_{x=S(u)}\frac{\partial S(u)}{\partial u}+\left.\frac{\partial f(x,u)}{\partial u}\right|_{x=S(u)}=0.
\end{equation}
Next, it follows from the definition (\ref{Tdef}) that
\begin{equation}\label{y}
T'(u)=\left.\frac{\partial h(x)}{\partial x}\right|_{x=S(u)}\left.\frac{\partial S(u)}{\partial u}.\right.
\end{equation}
Combining (\ref{x}) and (\ref{y}), and substituting (\ref{abcdef}), we verify (\ref{acdc}). Substituting (\ref{acdc}), we then rewrite (\ref{stepby}) as
\begin{equation}
-I+\lambda_kT'({\mathbf u}^*),
\end{equation}
and conclude that it indeed has a positive eigenvalue, because $\lambda_k<0$ implies that $\lambda_kT'({\mathbf u}^*)$ is a nonnegative matrix and (\ref{unstabk}) implies that its spectral radius exceeds one. Since the spectral radius is an eigenvalue for nonnegative matrices (see, {\it e.g.}, \cite{berman}), the conclusion follows.
\hfill $\Box$
\medskip

The eigenvectors $v_k$ of $P$ used in the similarity transformation (\ref{similarity}) may be interpreted as the spatial modes of the system. Thus, the stability properties of the matrix (\ref{eachk}) for each $k$ determines whether the corresponding mode decays or grows in time.  Since the spectral radius is nonnegative and $\lambda_k$, $k=1,\cdots,N$, are in decreasing order, whenever the instability criterion (\ref{unstabk}) holds for a particular mode $k$, it also holds for higher values of $k$. Because larger wavenumbers $k$ imply higher spatial frequency content in $v_k$, we conclude that the instability condition above sets the stage for the formation of high-frequency spatial patterns.
\smallskip

\section{Patterning in Bipartite Graphs}\label{BIP}
\subsection{Emergence of Checkerboard Patterns}
For bipartite graphs, where $\lambda_N=-1$ as stated in Lemma \ref{prop}, the instability condition in Theorem \ref{T1} is:
\begin{equation}\label{unstabbipartite}
\rho(T'({\mathbf u}^*))>1.
\end{equation}
This condition indicates the growth of the highest spatial-frequency mode $v_N$ which exhibits opposite signs for adjacent nodes. Thus, concentrations in adjacent nodes move in opposite directions in the vicinity of the homogeneous steady-state. We now show that, if the map
\begin{equation}\label{square}
T^2(\cdot):=T(T(\cdot))
\end{equation}
has two fixed points ${\mathbf u}_1 \neq {\mathbf u}_2$ other than ${\mathbf u}^*$, satisfying:
\begin{equation}\label{period2}
{\mathbf u}_1=T({\mathbf u}_2), \quad {\mathbf u}_2=T({\mathbf u}_1),
\end{equation}
then the system (\ref{one})-(\ref{two}) has an inhomogeneous steady-state with two sets of concentrations, each assigned to one of two adjacent cells.  We will refer to this steady-state as a
 ``checkerboard" pattern, since adjacent cells adopt distinct states. Although this term may be associated with cells arranged as a grid graph in two dimensional space, we will use it broadly for any spatial arrangement that forms a bipartite graph.

\begin{proposition}\label{sss}
Let $\mathcal{G}$ be a bipartite graph and let the sets
$\mathcal{I}\subset \{1,\cdots,N\}$ and $\mathcal{I}'=\{1,\cdots,N\}-\mathcal{I}$ be such that no two nodes in the same set are adjacent. If there exist ${\mathbf u}_1 \in \mathscr{U}$ and ${\mathbf u}_2 \in \mathscr{U}$, ${\mathbf u}_1\neq {\mathbf u}_2$, satisfying (\ref{period2}), then
\begin{equation}\label{onoff}
x^i=S({\mathbf u}_1), \ i \in \mathcal{I}, \quad  x^i=S({\mathbf u}_2), \ i \in \mathcal{I'},
\end{equation}
and
\begin{equation}\label{offon}
x^i=S({\mathbf u}_2), \ i \in \mathcal{I}, \quad  x^i=S({\mathbf u}_1), \ i \in \mathcal{I'},
\end{equation}
are steady-states for system (\ref{one})-(\ref{two}).
\end{proposition}

\noindent
{\it Proof:}
To show that (\ref{onoff}) is a steady-state, we note that, if $i\in \mathcal{I}$, then $y^i=T({\mathbf u}_1)$ and, if $i\in \mathcal{I}'$, then $y^i=T({\mathbf u}_2)$.
From (\ref{two}), the input $u^i$ to a node in  $\mathcal{I}$ is $T({\mathbf u}_2)$ because all neighbors of this node belong to $\mathcal{I}'$. Likewise, the input $u^i$ to a node in  $\mathcal{I}'$ is $T({\mathbf u}_1)$ because all neighbors of this node belong to $\mathcal{I}$. Since $T({\mathbf u}_2)={\mathbf u}_1$ and $T({\mathbf u}_1)={\mathbf u}_2$, we conclude  that (\ref{onoff}) is indeed a steady-state, and identical arguments apply to (\ref{offon}). \hfill $\Box$

\begin{theorem}\label{as}
Consider the system (\ref{one})-(\ref{two}) and suppose Assumptions \ref{characteristic} and \ref{iomon}, and the hypotheses of Proposition \ref{sss} hold. If, in addition,
 \begin{equation}\label{key}
\rho(T'({\mathbf u}_1)T'({\mathbf u}_2))<1,
\end{equation}
 then the steady-states (\ref{onoff}) and (\ref{offon}) are asymptotically stable.
\end{theorem}

Before giving the proof, we note that (\ref{period2}) corresponds to a period-two orbit of the discrete-time system:
\begin{equation}\label{dt}
u(t+1)=T(u(t)),
 \end{equation}
 and (\ref{key}) implies the asymptotic stability of this orbit.
 Likewise, (\ref{unstabbipartite}) indicates instability of the fixed point ${\mathbf u}^*$ for this discrete-time system. Thus, an interesting duality exists between (\ref{dt}) and the spatially-distributed system (\ref{one})-(\ref{two}) defined on a bipartite graph: A bifurcation from a stable fixed point to a stable period-two orbit in (\ref{dt}) corresponds to  the emergence of stable checkerboard patterns
from a homogeneous steady-state in (\ref{one})-(\ref{two}).

 In the single-input, single-output case with $\mathscr{U}=\mathbb{R}_{\ge 0}$, where $T:\mathbb{R}_{\ge 0}\rightarrow \mathbb{R}_{\ge 0}$ is a nonincreasing function by Assumption \ref{iomon}, condition (\ref{unstabbipartite}) indeed implies the existence of a period-two orbit (\ref{period2}).
 To see this, assume to the contrary that ${\mathbf u}^*$ is the unique fixed point of $T^2(\cdot)$.
Since $T(\cdot)$ is continuous and nonincreasing, this uniqueness property would imply that ${\mathbf u}^*$ is a global attractor for all solutions of the difference equation (\ref{dt}) starting in $\mathbb{R}_{\ge 0}$ \cite[Lemma 1.6.5]{kocic}. This, however, contradicts (\ref{unstabbipartite}), which implies instability of ${\mathbf u}^*$ for this scalar difference equation.

The argument above does not suggest the uniqueness of the pair $({\mathbf u}_1,{\mathbf u}_2)$, and multiple pairs satisfying (\ref{period2}) may exist. However,  we claim that at least one pair satisfies:
\begin{equation}\label{key2}
\left.\frac{dT^2(u)}{du}\right|_{u={\mathbf u}_1}=\left.\frac{dT^2(u)}{du}\right|_{u={\mathbf u}_2}=T'({\mathbf u}_1)T'({\mathbf u}_2)<1,
\end{equation}
which is the scalar equivalent of (\ref{key}), since $T'({\mathbf u}_1)T'({\mathbf u}_2)$ is nonnegative.
To see this, note from (\ref{unstabbipartite}) that:
\begin{equation}
\left.\frac{dT^2(u)}{du}\right|_{u={\mathbf u}^*}=T'({\mathbf u}^*)T'({\mathbf u}^*)>1
\end{equation}
and suppose, in contrast to (\ref{key2}), that the derivative of $T^2(\cdot)$ is greater than or equal to one at each of its fixed points. This implies that $T^2(u)\ge u$ for all $u\ge {\mathbf u}^*$, because $T^2(u)-u$ has nonnegative slope at zero-crossings and, thus, remains nonnegative for $u\ge {\mathbf u}^*$.
The inequality $T^2(u)\ge u$ implies unbounded growth of  $T^2(\cdot)$ which is a contradiction because $T(\cdot)$ is continuous and nonincreasing, thus, bounded.

\bigskip

\noindent
{\bf Proof of Theorem \ref{as}:} Let $N_\mathcal{I}$ and $N_\mathcal{I'}:=N-N_\mathcal{I}$ denote the cardinalities of the sets $\mathcal{I}$ and $\mathcal{I}'$, and index the cells such that $i=1,\cdots,N_\mathcal{I}$ belong to $\mathcal{I}$, and $i=N_\mathcal{I}+1,\cdots,N$ belong to $\mathcal{I}'$.
Then the matrix $P$ has the form:
\begin{equation}\label{specialP}
P=\left[\begin{array}{cc} 0 & P_{12} \\ P_{21} & 0 \end{array}\right]
\end{equation}
where $P_{12}\in \mathbb{R}^{N_\mathcal{I}\times N_\mathcal{I'}}$, $P_{21}\in \mathbb{R}^{N_\mathcal{I'}\times N_\mathcal{I}}$. Let $X:=[{x^1}^T \cdots {x^N}^T]^T$, and note that the linearization of (\ref{one})-(\ref{two}) about (\ref{onoff}) gives the Jacobian matrix:
\begin{equation}\label{composite}
\left[ \begin{array}{cc} I_{N_\mathcal{I}}\otimes A_1 & P_{12}\otimes (B_1C_2) \\
P_{21}\otimes (B_2C_1) & I_{N_\mathcal{I'}}\otimes A_2\end{array} \right]
\end{equation}
where
\begin{equation}\label{AE}
A_j:=\left.\frac{\partial f(x,u)}{\partial x}\right|_{(x,u)=(S({\mathbf u}_j),{\mathbf u}_j)}, \quad B_j:=\left.\frac{\partial f(x,u)}{\partial u}\right|_{(x,u)=(S({\mathbf u}_j),{\mathbf u}_j)}, \quad C_j:=\left.\frac{\partial h(x)}{\partial x}\right|_{x=S({\mathbf u}_j)}, \quad j=1,2.
\end{equation}
From the definition (\ref{randomwalk}),  the matrix $DP$, where $D$ is a diagonal matrix of the node degrees, is symmetric. Since $D^{-1/2}(DP)D^{-1/2}=D^{1/2}PD^{-1/2}$
is also symmetric, we write:
\begin{equation}\label{BB}
D^{1/2}PD^{-1/2}=\left[\begin{array}{cc} 0 & R \\ R^T & 0 \end{array}\right]
\end{equation}
where $R\in \mathbb{R}^{N_\mathcal{I}\times N_\mathcal{I'}}$ is appropriately defined. Then, we apply the following similarity transformation to (\ref{composite}):
\begin{equation}\label{compositetilde}
(D^{1/2}\otimes I_n) \left[ \begin{array}{cc} I_{N_\mathcal{I}}\otimes A_1 & P_{12}\otimes (B_1C_2) \\
P_{21}\otimes (B_2C_1) & I_{N_\mathcal{I'}}\otimes A_2\end{array} \right](D^{-1/2}\otimes I_n)=\left[ \begin{array}{cc} I_{N_\mathcal{I}}\otimes A_1 & R\otimes (B_1C_2) \\
R^T\otimes (B_2C_1) & I_{N_\mathcal{I'}}\otimes A_1\end{array} \right].
\end{equation}
The structure of (\ref{BB}) is such that it can diagonalized with an orthonormal matrix of the form:
\begin{equation}\label{specialV}
Q=\left[ \begin{array}{rrrr} Q_1 & Q_1 & Q_3 & 0 \\ Q_2 & -Q_2 & 0 & Q_4\end{array}\right]
\end{equation}
which results in:
\begin{equation}\label{diagd}
\left[\begin{array}{cc} 0 & R \\ R^T & 0 \end{array}\right]Q=Q\left[ \begin{array}{rrrr} \Lambda_+ &  & & \\ & -\Lambda_+ & & \\ & & 0 & \\ & & & 0\end{array} \right]
\end{equation}
where $\Lambda_+$ is a diagonal matrix of the strictly positive eigenvalues of $P$, the columns of $Q_3$ and $Q_4$ span the null spaces of $R^T$ and $R$, respectively,
and the dimensions of the zero diagonal blocks in (\ref{diagd}) are consistent with the dimensions of these null spaces (which we denote as $n_3$ and $n_4$, respectively).
From the orthonormality of $Q$, we get the identities:
\begin{eqnarray}\label{first}
&& Q_1^TQ_1=Q_2^TQ_2=\frac{1}{2}I_r\\
&&Q_4^TQ_4=I_{n_4} \quad Q_3^TQ_3=I_{n_3}\\
&&Q_1^TQ_3=0 \quad Q_2^TQ_4=0,\label{mid}
\end{eqnarray}
where $r$ is the dimension of $\Lambda_+$. Likewise, equation (\ref{diagd}) implies:
\begin{eqnarray}
 RQ_2=Q_1\Lambda_+ && R^TQ_1=Q_2\Lambda_+\\
RQ_4=0 && R^TQ_3=0. \label{last}
\end{eqnarray}
We now return to the Jacobian matrix (\ref{compositetilde}) and further apply the following similarity transformation:
\begin{equation}\label{produ}
\left[ \!\!\!
\begin{array}{cc} 2Q_1^T \otimes I_n & 0 \\
0 & 2Q_2^T \otimes I_n \\
Q_3^T \otimes I_n & 0 \\
0 & Q_4^T \otimes I_n
\end{array}
\!\!\!\right]
\left[ \begin{array}{cc} I_{N_\mathcal{I}}\otimes A_1 & R\otimes (B_1C_2) \\
R^T\otimes (B_2C_1) & I_{N_\mathcal{I'}}\otimes A_1\end{array} \right]
\left[\!\!\!
\begin{array}{cccc}
Q_1\otimes I_n & 0 & Q_3\otimes I_n & 0\\
0 & Q_2\otimes I_n & 0 & Q_4\otimes I_n
\end{array}
\!\!\!\right]
\end{equation}
where the leftmost matrix is the inverse of the rightmost matrix from (\ref{first})-(\ref{mid}). Likewise,
using (\ref{first})-(\ref{last}), it is not difficult to show that the product (\ref{produ}) equals:
\begin{equation}\label{almost}
\left[\begin{array}{cccc}
I_r\otimes A_1 & \Lambda_+\otimes (B_1C_2) & &\\
\Lambda_+\otimes (B_2C_1) & I_r \otimes A_{2} & & \\
& &  I_{n_3}\otimes A_{1} & \\ & & & I_{n_4}\otimes A_{2}
\end{array}\right].
\end{equation}
Since Assumptions \ref{characteristic} and \ref{iomon} imply that $A_1$ and $A_2$ are Hurwitz \cite[Lemma 6.5]{AngSon04}, stability of (\ref{almost}) is determined by the upper left blocks which, upon a similarity transformation with an appropriate permutation matrix, are block-diagonalized into $r$ blocks of the form:
\begin{equation}\label{blok}
\left[ \begin{array}{cc} A_1 & \lambda_i B_1C_2 \\ \lambda_i B_2C_1 & A_2 \end{array}\right]
\end{equation}
$i=1,\cdots,r$.

We will now show that (\ref{blok}) is Hurwitz for any $\lambda_i\in [-1,1]$.  Since all eigenvalues of $P$ lie in this interval by Lemma \ref{prop}, this will conclude the proof.  We do not provide a separate proof for the asymptotic stability of (\ref{offon}), as identical arguments apply when the indices $1$ and $2$ are swapped in (\ref{blok}). If $\lambda_i=0$, (\ref{blok}) is Hurwitz because $A_1$ and $A_2$ are Hurwitz. If  $\lambda_i\neq 0$, then we apply the similarity transformation:
\begin{equation}
\left[\begin{array}{cc} I & 0 \\ 0 & \lambda_i^{-1}I\end{array} \right]\left[ \begin{array}{cc} A_1 & \lambda_i B_1C_2 \\ \lambda_i B_2C_1 & A_2 \end{array}\right]\left[\begin{array}{cc} I & 0 \\ 0 & \lambda_iI\end{array} \right]=\left[ \begin{array}{cc} A_1 & \lambda_i^2 B_1C_2 \\ B_2C_1 & A_2 \end{array}\right]
\end{equation}
and rewrite the result as:
\begin{equation}\label{ABC}
\mathcal{A}+\mathcal{BC}
\end{equation}
where
\begin{equation}\label{defABC}
\mathcal{A}:=\left[ \begin{array}{cc} A_1 & \lambda_i^2 B_1C_2 \\ 0 & A_2 \end{array}\right], \quad \mathcal{B}=\left[ \begin{array}{c} 0 \\ B_2 \end{array}\right],\quad \mathcal{C}=\left[\ C_1\ \ 0 \ \right].
\end{equation}
We claim that the linear system defined by the triplet $(\mathcal{C},\mathcal{A},\mathcal{B})$ is monotone with respect to $K^U=K^Y=\mathbb{R}^m_{\ge 0}$, and $K^X=-K\times K$ where $K$ is as in Assumption \ref{iomon}. To see this, first note from Lemma \ref{jaclin} that $(C_1,A_1,B_1)$ and $(C_2,A_2,B_2)$ are monotone with respect to the cones specified in Assumption \ref{iomon}. By Lemma \ref{linchar}, this means that:
\begin{equation}\label{useit}
x\in K \ \Rightarrow \ A_jx \in K, \quad u\in \mathbb{R}^m_{\ge 0} \ \Rightarrow \ B_j\,u \in K, \quad x\in K \ \Rightarrow \ C_jx \in \mathbb{R}^m_{\le 0}, \quad j=1,2.
\end{equation}
We now show that the conditions of Lemma \ref{linchar} hold for $(\mathcal{C},\mathcal{A},\mathcal{B})$ with $K^U=K^Y=\mathbb{R}^m_{\ge 0}$, $K^X=-K\times K$:

\noindent
1) Suppose $x=[x_1^T x_2^T]^T\in -K\times K$, that is $x_1\in -K$, $x_2\in K$. Then,
\begin{equation}
\mathcal{A}x=\left[\begin{array}{cc}A_1x_1+\lambda_i^2B_1C_2x_2 \\ A_2x_2\end{array}\right]\in -K\times K
\end{equation}
because, from (\ref{useit}), $A_1x_1\in -K$, $A_2x_2\in K$, $C_2x_2\in \mathbb{R}^m_{\le 0}$ and, hence, $B_1C_2x_2\in -K$.

\noindent
2) We want to show that $u\in \mathbb{R}^m_{\ge 0}$ implies $\mathcal{B}u \in -K\times K$. From the definition of $\mathcal{B}$ in (\ref{defABC}), $\mathcal{B}u \in -K\times K$ means  $B_2u\in K$. It follows from the second implication in (\ref{useit}) that $u\in \mathbb{R}^m_{\ge 0}$ indeed implies $B_2u\in K$.

\noindent
3) To prove monotonicity with $K^Y=\mathbb{R}^m_{\ge 0}$, we need to show that $x_1\in -K$ and $x_2\in K$ imply $\mathcal{C}\,[x_1^T x_2^T]^T\in \mathbb{R}^m_{\ge 0}$. This is indeed true, since $\mathcal{C}\,[x_1^T x_2^T]^T=C_1x_1$ and, from (\ref{useit}), $x_1\in -K$ implies $C_1x_1\in \mathbb{R}^m_{\ge 0}$.

Having verified the conditions of Lemma \ref{linchar}, we conclude that $(\mathcal{C},\mathcal{A},\mathcal{B})$ is monotone with respect to $K^U=K^Y=\mathbb{R}^m_{\ge 0}$. In addition, the matrix $\mathcal{A}$ in (\ref{defABC}) is  Hurwitz, as $A_1$ and $A_2$ are Hurwitz. Thus, it follows from the first statement in Lemma \ref{posfbk} that, if $-(I+\mathcal{C}\mathcal{A}^{-1}\mathcal{B})$ is Hurwitz, then so is (\ref{ABC}).
Note that
\begin{equation}\label{alm}
\mathcal{C}\mathcal{A}^{-1}\mathcal{B}=\left[\ C_1\ \ 0 \ \right]\left[ \begin{array}{cc} A_1^{-1} & -\lambda_i^2 A_1^{-1}B_1C_2A_2^{-1} \\ 0 & A_2^{-1} \end{array}\right]\left[ \begin{array}{c} 0 \\ B_2 \end{array}\right]=-\lambda_i^2 C_1A_1^{-1}B_1C_2A_2^{-1}B_2
\end{equation}
and, from a derivation similar to the one for (\ref{acdc}), $T'({\mathbf u}_j)=-C_jA_j^{-1}B_j$, $j=1,2.$ Thus, (\ref{alm}) gives:
 \begin{equation}
-(I+\mathcal{C}\mathcal{A}^{-1}\mathcal{B})=-I+\lambda_i^2T'({\mathbf u}_1)T'({\mathbf u}_2),
\end{equation}
and (\ref{key}) and $\lambda_i\in [-1,1]$ imply that $-(I+\mathcal{C}\mathcal{A}^{-1}\mathcal{B})$ is indeed Hurwitz. From Lemma \ref{posfbk}, this means that (\ref{ABC}) and, thus, (\ref{blok}) is Hurwitz $i=1,\cdots,r,$ concluding the proof.
\mbox{} \hfill $\Box$

\subsection{Generic Convergence to Steady-States}
Thus far we have studied {\it local} asymptotic stability properties of the steady-states.  Strongly monotone systems (as defined in the paragraph above Definition \ref{angson}) have been shown to possess a ``generic convergence" property \cite{Hirsch-Smith,smith} which means that almost every bounded solution (except for a measure-zero set of initial conditions) converges to the set of steady-states.
Below we first prove monotonicity of (\ref{one})-(\ref{two}) in Theorem \ref{mono} and, next establish strong monotonicity in Theorem \ref{smon}, thereby concluding generic convergence for this system.

\begin{figure}[h]
\vspace{-1.8cm}
\begin{center}
\mbox{}\setlength{\unitlength}{1.3mm}
\begin{picture}(60,70)
\put(-20,20){\psfig{figure=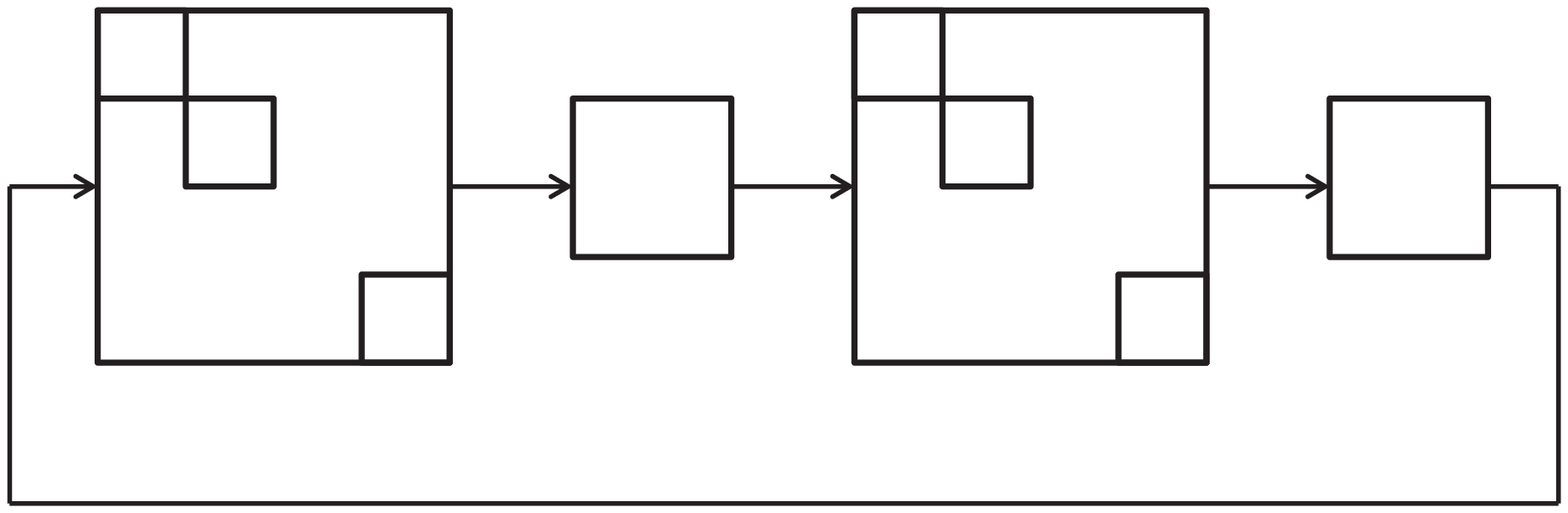,width=100\unitlength}}
 \put(-7.5,48.5){$x^1$}\put(-2.5,43.5){$x^2$}\put(2.5,38.5){$\ddots$}\put(7.0,33.5){$x^{N_\mathcal{I}}$}
 \put(18.5,41.5){$P_{21}\otimes I$}
  \put(35.5,48.5){}\put(40.5,43.5){}\put(45.5,38.5){$\ddots$}\put(50,33.5){$x^{N}$}
 \put(61.5,41.5){$P_{12}\otimes I$}
 \put(-20.5,43){$\mathcal{U}:=U^\mathcal{I}$}\put(13,43){$Y^\mathcal{I}$}
 \put(28.5,43){$U^\mathcal{I'}$}\put(56,43){$Y^\mathcal{I'}$}
\put(75,43){$\mathcal{Y}$}
\end{picture}
\vspace{-3.2cm}
 \caption{\small A block diagram for the system (\ref{one})-(\ref{two}) when the contact graph is bipartite and the corresponding interconnection matrix $P$ is decomposed as in (\ref{specialP}).}  \label{bipardecomp}
\end{center}
\vspace{-.5cm}
\end{figure}

\begin{theorem}\label{mono} If $\mathcal{G}$ is bipartite and Assumption  \ref{iomon} holds, then the system (\ref{one})-(\ref{two}) is monotone.
\end{theorem}
\noindent
{\it Proof:}
Let $\mathcal{I}\subset \{1,\cdots,N\}$ and $\mathcal{I}'=\{1,\cdots,N\}-\mathcal{I}$ be defined as in Proposition \ref{sss}, and suppose that in (\ref{one}), the cells are indexed such that $i=1,\cdots,N_\mathcal{I}$ belong to $\mathcal{I}$, and $i=N_\mathcal{I}+1,\cdots,N$ belong to $\mathcal{I}'$ as in the proof of Theorem \ref{as}, where $N_\mathcal{I}$ is the cardinality of set $\mathcal{I}$. Let $X^\mathcal{I}:=[{x^1}^T \cdots {x^{N_\mathcal{I}}}^T]^T$, $X^\mathcal{I'}:=[{x^{N_\mathcal{I}+1}}^T \cdots {x^N}^T]^T$,  and define $U^\mathcal{I}$, $U^\mathcal{I'}$, $Y^\mathcal{I}$, $Y^\mathcal{I'}$ similarly. Then, the interconnection condition (\ref{two}) becomes:
\begin{eqnarray}\label{twoa}
U^\mathcal{I}&=&(P_{12}\otimes I_m)Y^\mathcal{I'} \\
\label{twob}
 U^\mathcal{I'}&=&(P_{21}\otimes I_m)Y^\mathcal{I}
\end{eqnarray}
where $P_{12}$ and $P_{21}$ are as in (\ref{specialP}). A block diagram illustrating this interconnection is depicted in Figure \ref{bipardecomp}.

To prove the monotonicity of this feedback system, we establish the monotonicity of the feedforward system with input $\mathcal{U}:=U^\mathcal{I}$ and output $\mathcal{Y}:=(P_{12}\otimes I_m)Y^\mathcal{I'}$:
\smallskip

\noindent
{\it Claim: The feedforward system in Figure \ref{bipardecomp} with input $\mathcal{U}$ and output $\mathcal{Y}$ is monotone with respect to the positivity cones
$K^U=K^Y=\mathbb{R}^{mN_\mathcal{I}}_{\ge 0}$, and $K^X=K^{N_\mathcal{I}}\times \{-K\}^{N-N_\mathcal{I}}$.}
\smallskip

\noindent
The theorem follows from this claim because a monotone input-output system, where the inputs and outputs are ordered with respect to the same positivity cone, is monotone when the output is connected to the input with unitary positive feedback $\mathcal{U}=\mathcal{Y}$ (see the first part of the proof of \cite[Theorem 2]{AngSon04}).

To prove the claim above, we take two input signals satisfying $\mathcal{U}(t)\preceq \hat{\mathcal{U}}(t)$ for all $t\ge 0$, which means that $u^i(t)\preceq \hat{u}^i(t)$, $i\in \mathcal{I}$, with respect to $\mathbb{R}^m_{\ge 0}$.
 Likewise, we let $X(0)\preceq \hat{X}(0)$ with respect to the cone $K^X=K^{N_\mathcal{I}}\times \{-K\}^{N-N_\mathcal{I}}$, which means that $x^i(0)\preceq \hat{x}^i(0)$ for $i\in \mathcal{I}$ and $x^i(0)\succeq \hat{x}^i(0)$ for $i\in \mathcal{I'}$ with respect to the cone $K$.
It follows from Assumption \ref{iomon} that:
\begin{equation}\label{piece1}
x^i(t)\preceq \hat{x}^i(t) \quad \forall t\ge 0 \quad  i\in \mathcal{I}.
\end{equation}
Moreover, since $x\preceq \hat{x}$ implies $h(x)\preceq h(\hat{x})$ with respect to $\mathbb{R}^m_{\le 0}$ by Assumption \ref{iomon}, we conclude ${Y}^\mathcal{I}(t)\succeq \hat{Y}^\mathcal{I}(t)$ with respect to $\mathbb{R}_{\ge 0} ^{mN_\mathcal{I}}$. Because $P_{21}$ is a nonnegative matrix, (\ref{twob}) implies ${U}^\mathcal{I'}(t)\succeq \hat{U}^\mathcal{I'}(t)$ which means that $u^i(t)\succeq \hat{u}^i(t)$ for all $t\ge 0$, $i\in \mathcal{I'}$. As noted above, $x^i(0)\succeq \hat{x}^i(0)$ for $i\in \mathcal{I'}$ and, hence, another application of Assumption \ref{iomon} yields:
 \begin{equation}\label{piece2}
x^i(t)\succeq \hat{x}^i(t) \quad \forall t\ge 0 \quad  i\in \mathcal{I'}.
\end{equation}
Since (\ref{piece1}) and (\ref{piece2}) hold with respect to $K$, we conclude that $X(t)\preceq \hat{X}(t)$ for all $t\ge 0$ with respect to $K^X=K^{N_\mathcal{I}}\times \{-K\}^{N-N_\mathcal{I}}$.
To conclude the proof of the claim, we need to show that $X \preceq \hat{X}$ implies $\mathcal{Y} \preceq \hat{\mathcal{Y}}$. Indeed, the former implies $x^i \succeq \hat{x}^i$ for $i\in \mathcal{I'}$ and, it follows from Assumption \ref{iomon} that $h(x^i) \succeq h(\hat{x}^i)$ with respect to $R^m_{\le 0}$. Thus, ${Y}^\mathcal{I'} \preceq \hat{Y}^\mathcal{I'}$ with respect to $\mathbb{R}_{\ge 0} ^{m(N-N_\mathcal{I})}$ and, since $P_{12}$ is a nonnegative matrix, we conclude $\mathcal{Y} \preceq \hat{\mathcal{Y}}$ with respect to $\mathbb{R}_{\ge 0} ^{mN_\mathcal{I}}$. \hfill $\Box$

\bigskip

To establish {\it strong} monotonicity, we need additional {\it excitability} and {\it transparency} conditions,
as defined in \cite{AngSon04,EncSon05}:

\begin{defn}\label{extran}
The monotone system  $\dot{x}=f(x,u)$, $y=h(x)$  is said to be {\it excitable} if $x(0)\preceq \hat{x}(0)$ and $u(t)\prec \hat{u}(t)$ for almost all $t>0$ imply $x(t)\ll \hat{x}(t)$ $\forall t>0$. It is said to be {\it transparent} if  $u(t)\preceq \hat{u}(t)$ and $x(0)\prec \hat{x}(0)$ imply $y(t)\ll \hat{y}(t)$ $\forall t>0$.
\end{defn}

Since inputs and outputs are ordered with respect to orthants ($K^U=\mathbb{R}^m_{\ge 0}$ and $K^Y=-K^U$) in Assumption \ref{iomon}, here we give a less restrictive definition of excitability (transparency) which requires that this property hold with respect to a particular component of the input (output) vector:

\begin{defn}\label{extran-new}
The monotone system  $\dot{x}=f(x,u)$, $y=h(x)$  is said to be {\it excitable} by the $k$th input if $x(0)\preceq \hat{x}(0)$, $u(t)\preceq \hat{u}(t)$ and $u_k(t)\prec \hat{u}_k(t)$ for almost all $t>0$ imply $x(t)\ll \hat{x}(t)$ $t>0$. It is said to be {\it transparent} from the $k$th output if $u(t)\preceq \hat{u}(t)$ and $x(0)\prec \hat{x}(0)$ imply $y_k(t)\prec \hat{y}_k(t)$ $\forall t>0$.
\end{defn}

\begin{assumption}\label{either}
There exists $k\in \{1,\cdots,m\}$ such that system (\ref{one}) is excitable by the $k$th input and transparent from the $k$th output.
\end{assumption}

\begin{theorem}\label{smon}
If, in addition to the conditions of Theorem \ref{mono}, Assumption \ref{either} holds, then (\ref{one})-(\ref{two}) is strongly monotone.
\end{theorem}

\noindent
{\it Proof:}  We need to show that $X(0)\prec \hat{X}(0)$ implies $X(t)\ll \hat{X}(t)$ for all $t>0$ with respect to the cone $K^X=K^{N_\mathcal{I}}\times \{-K\}^{N-N_\mathcal{I}}$ for which monotonicity was proven in Theorem \ref{mono}.  By this monotonicity property, we already know that $X(t)\preceq \hat{X}(t)$ for all $t\ge 0$, and Assumption \ref{iomon} implies $y^i(t) \succeq \hat{y}^i(t)$ if $i\in \mathcal{I}$, $y^i(t) \preceq \hat{y}^i(t)$ if $i\in \mathcal{I'}$, both with respect to $\mathbb{R}^m_{\ge 0}$. Because $P_{12}$ and $P_{21}$ in (\ref{twoa})-(\ref{twob}) are nonnegative matrices, we conclude:
\begin{equation}\label{ref1}
u^i(t) \preceq \hat{u}^i(t) \ \, i\in \mathcal{I}, \quad u^i(t) \succeq \hat{u}^i(t) \ \, i\in \mathcal{I'}.
\end{equation}
Next, note that $X(0)\prec \hat{X}(0)$ means $x^i(0)\neq \hat{x}^i(0)$ for at least one $i\in \{1,\cdots,N\}$, say $i^*$. Thus, with respect to the cone $K$:
\begin{equation}\label{ref2}
x^{i^*}(0)\prec \hat{x}^{i^*}(0) \ \, \mbox{if} \, \ i^*\in \mathcal{I}, \quad x^{i^*}(0)\succ \hat{x}^{i^*}(0) \ \, \mbox{if} \ \, i^*\in \mathcal{I'}.
\end{equation}
Using (\ref{ref1})-(\ref{ref2}) and the transparency assumption from the $k$th output, we conclude that the following holds for all $t>0$ with respect to the standard order induced by $\mathbb{R}_{\ge 0}$:
\begin{equation}\label{ref3}
y_k^{i^*}(t) \succ \hat{y}_k^{i^*}(t) \ \, \mbox{if} \, \ i^*\in \mathcal{I}, \quad y_k^{i^*}(t) \prec \hat{y}_k^{i^*}(t) \ \, \mbox{if} \, \ i^*\in \mathcal{I'}.
\end{equation}
Now, pick an arbitrary $i^\dagger \in \{1,\cdots,N\}$, and note from the connectedness of the contact graph $\mathcal{G}$ that a path of adjacent cells $i_\ell$, $\ell=1,\cdots,M$ exists such that  $i_1=i^*$ and $i_M=i^\dagger$. Since $i_2$ is a neighbor of $i_1=i^*$, for all $t>0$,
\begin{equation}\label{ref4}
u_k^{i_2}(t) \succ \hat{u}_k^{i_2}(t) \ \, \mbox{if} \, \ i^*\in \mathcal{I}, \quad u_k^{i_2}(t) \prec \hat{u}_k^{i_2}(t) \ \, \mbox{if} \, \ i^*\in \mathcal{I'}.
\end{equation}
Since $\mathcal{G}$ is bipartite, $i^*\in \mathcal{I}$ means $i_2\in \mathcal{I'}$, and $i^*\in \mathcal{I'}$ means $i_2\in \mathcal{I}$. Thus, from $X(0)\prec \hat{X}(0)$:
\begin{equation}\label{ref5}
x^{i_2}(0)\succeq \hat{x}^{i_2}(0) \ \, \mbox{if} \, \ i^*\in \mathcal{I}, \quad x^{i_2}(0)\preceq \hat{x}^{i_2}(0) \ \, \mbox{if} \ \, i^*\in \mathcal{I'}.
\end{equation}
From the excitability assumption by the $k$th input, (\ref{ref4}) and (\ref{ref5}) imply:
\begin{equation}\label{ref6}
x^{i_2}(t)\gg \hat{x}^{i_2}(t) \ \, \mbox{if} \, \ i^*\in \mathcal{I}, \quad x^{i_2}(t)\ll \hat{x}^{i_2}(t) \ \, \mbox{if} \ \, i^*\in \mathcal{I'}
\end{equation}
$\forall t>0$  and, from transparency, the following holds with respect to the standard order:
\begin{equation}\label{ref7}
y_k^{i_2}(t)\prec \hat{y}_k^{i_2}(t) \ \, \mbox{if} \, \ i^*\in \mathcal{I}, \quad y_k^{i_2}(t)\succ \hat{y}_k^{i_2}(t) \ \, \mbox{if} \ \, i^*\in \mathcal{I'}.
\end{equation}
Continuing recursively, we conclude that (\ref{ref4})-(\ref{ref7}) hold for $i_\ell$, $\ell=3,\cdots,M$, with the inequalities reversed when $\ell$ is odd. In particular,  (\ref{ref6}) becomes:
\begin{equation}\label{ref8}
(-1)^\ell x^{i_\ell}(t)\gg (-1)^\ell\hat{x}^{i_\ell}(t) \ \, \mbox{if} \, \ i^*\in \mathcal{I}, \quad (-1)^\ell x^{i_\ell}(t)\ll (-1)^\ell\hat{x}^{i_\ell}(t) \ \, \mbox{if} \ \, i^*\in \mathcal{I'}.
\end{equation}
Since $\mathcal{G}$ is bipartite,
if $M$ is even, $i^*=i_1\in \mathcal{I}$ means $i^\dagger = i_M\in \mathcal{I'}$, and $i^*\in \mathcal{I'}$ means $i^\dagger \in \mathcal{I}$. Likewise, if $M$ is odd, $i^*\in \mathcal{I}$ means $i^\dagger \in \mathcal{I}$, and $i^*\in \mathcal{I'}$ means $i^\dagger \in \mathcal{I'}$. Thus, (\ref{ref8}) with $\ell=M$ gives:
\begin{equation}\label{ref9}
x^{i^\dagger}(t)\gg \hat{x}^{i^\dagger}(t) \ \, \mbox{if} \, \ i^\dagger\in \mathcal{I'}, \quad x^{i^\dagger}(t)\ll \hat{x}^{i^\dagger}(t) \ \, \mbox{if} \ \, i^\dagger\in \mathcal{I}.
\end{equation}
Since this inequality holds for each $i^\dagger \in \{1,\cdots,N\}$, we conclude $X(t)\ll \hat{X}(t)$ as desired. \hfill $\Box$
\medskip


In preparation for the examples in the next section, we now review a graphical test to ascertain excitability and transparency, given in \cite{AngSon04}. Suppose the system $\dot{x}=f(x,u)$, $y=h(x)$, $x\in \mathscr{X}\subset \mathbb{R}^n$,  $u\in \mathscr{U}\subset \mathbb{R}^m$, $y\in \mathscr{Y}\subset \mathbb{R}^p$, is such that, for each $j\neq k$, $\partial f_j(x,u)/\partial x_k$ is either identically zero, strictly positive, or strictly negative for all $(x,u)\in \mathscr{X}\times \mathscr{U}$. Likewise,  $\partial f_j(x,u)/\partial u_k$ and $\partial h_k(x,u)/\partial x_j$ have the same {\it sign definiteness} property for each $j$ and $k$. Associate to this system a directed {\it incidence graph} with vertices $x_1,\cdots,x_n$, $u_1,\cdots,u_m$, $y_1,\cdots,y_p$. A directed edge is drawn from $x_k$ to $x_j$, $j\neq k$, if $\partial f_j(x,u)/\partial x_k$ is nonzero, from $u_k$ to $x_j$ if $\partial f_j(x,u)/\partial u_k$ nonzero, and from $x_j$ to $y_k$ if $\partial h_k(x,u)/\partial x_j$ is nonzero. The following lemma, adapted\footnote{Theorems 4 and 5 in \cite{AngSon04} give analogous tests for excitability and transparency with respect to Definition \ref{extran}. Theorem 4 requires that each state be reachable from each input through a directed path, and Theorem 5 stipulates that a directed path exist from each state to each output. The statement in Lemma \ref{copy} for transparency from the $k$th output follows directly from Theorem 5, by taking $y_k$ to be the only output of the system.  The statement for excitability by the $k$th input follows from a straightforward modification of Theorem 4:  Read the second part of the proof of Theorem 4 by replacing $j^\star$ with $k$.} from \cite{AngSon04}, proves excitability and transparency for systems that are monotone with respect an orthant cone:

\begin{lemma}\label{copy} Suppose the system $\dot{x}=f(x,u)$, $y=h(x)$ is monotone with respect to an orthant cone and admits an incidence graph according to the rules described above. This system is excitable by the $k$th input if each state $x_j$ is reachable through a directed path from  $u_k$, and transparent from the $k$th output if a directed path exists from each state $x_j$ to  $y_k$.
\end{lemma}

\section{Examples}
\subsection{A Class of Systems that Encompasses the Notch Signaling Model of \cite{ColMonMaiLew96}}
As a special case of (\ref{one}), consider the single-input, single-output system:
\begin{eqnarray}
\dot{x}_1^i&=&-\gamma_1x_1^i+g_1(x_2^i) \nonumber \\
\nonumber
& \vdots &\\
\label{specialcase}
\dot{x}_j^i&=&-\gamma_jx_j^i+g_j(x_{j+1}^i)\\
\nonumber
& \vdots &\\
\dot{x}_n^i&=&-\gamma_nx_n^i+g_n(u^i) \nonumber\\
y^i&=&x_1^i \nonumber
\end{eqnarray}
where, for $j=1,\cdots,n$, $x^i_j\ge 0$ denotes the concentration of species $j$ in cell $i$, $\gamma_j>0$  represents the corresponding degradation rate, and $g_j:\mathbb{R}_{\ge 0}\rightarrow \mathbb{R}_{\ge 0}$ is a continuously differentiable function.

The reference \cite{ColMonMaiLew96} studied (\ref{specialcase}) for $n=2$ species, as a rough model for Notch signaling where the membrane-bound Delta ligands bind the Notch receptors in adjacent cells.  This leads to the cleavage of Notch and the release of its intracellular domain which then serves as a co-transcription factor that inhibits the production of Delta in the same cell. Thus, in (\ref{specialcase}), $x_1$ represents the concentration of Delta and $x_2$ represents the concentration of the co-transcription factor obtained from Notch. The function $g_1(\cdot)$ is assumed to be decreasing since the co-transcription factor inhibits the production of Delta, and $g_2(\cdot)$ is assumed to be increasing since Delta activates the production of the co-transcription factor in adjacent cells.

The reference \cite{ColMonMaiLew96} proved the emergence and stability of patterns for the case of $N=2$ cells, and observed the patterning behavior for $N>2$ by numerical simulations. A detailed bifurcation analysis is performed for this model in \cite{WeaShe01}, again for $N=2$.
In Proposition \ref{ACC} below, we show that the results of the present paper are applicable to the model (\ref{specialcase}) without restrictions on the number of species and cells.
In particular, the instability criterion for the homogeneous steady-state in Theorem \ref{T1}
makes use of the spectral properties of random walks and, unlike \cite{ColMonMaiLew96,Pla01} which analyze this steady-state for specific arrays, is applicable to arbitrary graphs.  Likewise, our study of checkerboard patterns in Proposition \ref{sss} and Theorem \ref{as} generalizes the statements in \cite{ColMonMaiLew96} for two cells to bipartite graphs of arbitrary size.
In addition, we establish monotonicity properties for bipartite graphs, thus revealing the global behavior of the solutions.

\begin{proposition}\label{ACC}
System (\ref{specialcase}) satisfies Assumption \ref{characteristic}. If an odd number of the functions $g_j(\cdot)$, $j=1,\cdots,n,$ are nonincreasing and the rest are nondecreasing\footnote{If one of the functions is constant, then one can count it as either nonincreasing or nondecreasing. However, this situation is of no interest in this paper, since the input-output characteristic (\ref{Tspecial}) is constant and, therefore, Theorems \ref{T1} and \ref{as} are not applicable.}, then it satisfies Assumption \ref{iomon} as well. If, in addition, the nondecreasing (nonincreasing) property is strengthened as:
\begin{equation}\label{strincdec}
g'_j(s)> 0 \quad (g'_j(s)< 0) \quad \forall s\ge 0 \quad j=1,\cdots,n,
\end{equation}
then Assumption \ref{either} also holds.
\end{proposition}

\noindent
{\it Proof:} We first prove that Assumption \ref{characteristic} holds. Given $u^*\ge 0$, the unique steady-state $x^*$ of (\ref{specialcase}) is given by:
\begin{equation}
x^*_n=\gamma_n^{-1}g_n(u^*), \quad x^*_j=\gamma_j^{-1}g_j(x^*_{j+1}), \ j=n-1,\cdots,1.
\end{equation}
In particular, the input-output characteristic is:
\begin{equation}\label{Tspecial}
T(\cdot):=\gamma_1^{-1}g_1(\gamma_2^{-1}g_2(\cdots(\gamma_n^{-1}g_n(\cdot)))).
\end{equation}
The Jacobian matrix:
\begin{equation}
\frac{\partial f(x,u)}{\partial x}=\left[ \begin{array}{cccc}-\gamma_1 & g'_1(x_2) & &  \\  & -\gamma_2 & \ddots & \\ &  & \ddots & g'_{n-1}(x_n)\\  & &  & -\gamma_n\end{array}\right]
\end{equation}
is upper-triangular with negative diagonal entries $-\gamma_j$, $j=1,\cdots,n$, and, hence, Hurwitz. This means that the determinant condition (\ref{det}) holds and the steady-state $x^*$ is asymptotically stable. Note from (\ref{specialcase}) that $x_n(t)$ exists for all $t\ge 0$ and converges to $\gamma_n^{-1}g_n(u^*)$. Applying a similar argument recursively for $j=n-1,\cdots,1$, we conclude that $x^*$ is {\it globally} asymptotically stable.

 To show that Assumption \ref{iomon} holds, we first select numbers $\epsilon_j\in \{0,1\}$, $j=1,\cdots,n$,  according to the following rule: Set $\epsilon_n=0$ if $g_n(\cdot)$ is nondecreasing, and $\epsilon_n=1$ if $g_n(\cdot)$ is nonincreasing. Then, for $j=n-1,n-2,\cdots,1$, set $\epsilon_j=\epsilon_{j+1}$ if $g_j(\cdot)$ is nondecreasing, and $\epsilon_j\neq \epsilon_{j+1}$ if $g_j(\cdot)$ is nonincreasing. It follows from this construction that, $\forall s\ge 0,$
\begin{equation}
(-1)^{\epsilon_n}g'_n(s)\ge 0, \qquad (-1)^{\epsilon_j+\epsilon_{j+1}}g'_j(s) \ge 0, \ j=1,\cdots,n-1.
\end{equation}
Since an odd number of the functions $g_j(\cdot)$ are nonincreasing, the selection of the numbers $\epsilon_j$ above yields $\epsilon_1=1$. Thus, an application of Lemma \ref{monotonetest} with $\delta_1=0$ and $\mu_1=1$ shows that the system (\ref{specialcase}) is monotone with respect to $K^U=\mathbb{R}_{\ge 0}$, $K^X=\{x\in \mathbb{R}^n \ | \ (-1)^{\epsilon_j}x_j\ge 0\}$, $K^Y=\mathbb{R}_{\le 0}$, as in Assumption \ref{iomon}.

To show that Assumption \ref{either} holds, we apply the test in Lemma \ref{copy}. The incidence graph for system (\ref{specialcase}) consists of the single path $u \mapsto x_n \mapsto x_{n-1} \mapsto \cdots \mapsto x_1 \mapsto y$, which means that any state is reachable from the input, and the output is reachable from any state. Thus, the system (\ref{specialcase}) is excitable and transparent.
\hfill $\Box$

\subsection{A Multi-Input, Multi-Output Model for Notch Signaling}
We now study the following system adapted\footnote{The equation corresponding to (\ref{cofactor}) in \cite{SprLakLeB11} includes a Hill function of $N^i\langle D^j \rangle_i$ instead of the linear term used here.} from the lateral inhibition model in \cite{SprLakLeB11}:
\begin{eqnarray}\label{receptor}
\dot{N}^i&=&\beta-\gamma N^i-k{N^i\langle D^j \rangle_i}\\
\dot{D}^i&=&g(S^i)-\gamma D^i-k{D^i\langle N^j \rangle_i}\\
\dot{S}^i&=&-\gamma S^i+k{N^i\langle D^j \rangle_i}.\label{cofactor}
\end{eqnarray}
Here, $N^i\ge 0$, $D^i\ge 0$, $S^i\ge 0$ are the concentrations in cell $i$ of the Notch receptor, Delta ligand, and a signaling protein activated by the binding on Delta and Notch, $k>0$, $\gamma>0$, $\beta>0$, $g:\mathbb{R}_{\ge 0} \rightarrow \mathbb{R}_{>0}$ is continuously differentiable, and decreasing since the production of Delta is inhibited by the signaling protein.  The notation
$\langle \cdot \rangle_i$ denotes the average of the quantity within brackets over all cells adjacent to  $i$. Unlike the model of \cite{ColMonMaiLew96} discussed in the previous subsection, (\ref{receptor}) incorporates the Notch receptor.

We let:
\begin{equation}
u_1^i:=\langle D^j \rangle_i, \ u_2^i:=\langle N^j \rangle_i, \ x_1^i=N^i, \ x_2^i=D^i, \ x_3^i=N^i+S^i,
\end{equation}
and rewrite (\ref{receptor})-(\ref{cofactor}) as:
\begin{eqnarray}\label{l1}
\dot{x}_1^i&=&\beta-\gamma x_1^i-k{x_1^iu_1^i}\\ \label{l2}
\dot{x}_2^i&=&g(x_3^i-x_1^i)-\gamma x_2^i-k{x_2^iu^i_2}\\ \label{l3}
\dot{x}_3^i&=&-\gamma x_3^i+\beta\\ \label{l4}
y_1^i&=&x_2^i\\ \label{l5}
y_2^i&=&x_1^i,
\end{eqnarray}
which is of the form (\ref{one}) with $\mathscr{X}=\{x\in \mathbb{R}^3 \ | \ x_1\ge 0, x_2\ge 0, x_3\ge x_1\}$, $\mathscr{U}=\mathscr{Y}=\mathbb{R}_{\ge 0}^2$.

\begin{proposition}
The system (\ref{l1})-(\ref{l5}), where $k>0$, $\gamma>0$, $\beta>0$, and $g:\mathbb{R}_{\ge 0} \rightarrow \mathbb{R}_{>0}$ is continuously differentiable, satisfies Assumption \ref{characteristic}.
If $g(\cdot)$ is nonincreasing, then it also satisfies Assumption \ref{iomon}. If $g'(s)>0$ for all $s\ge 0$, then Assumption \ref{either} holds for solutions in the forward invariant subset of $\mathscr{X}$ where $x_1>0$, $x_2 > 0$, $x_3=x_3^*$.
\end{proposition}

\noindent
{\it Proof:} Given $u_1^*\ge 0$, $u_2^*\ge 0$, the unique steady-state of (\ref{l1})-(\ref{l3}) is given by:
\begin{equation}
x_1^*=\frac{\beta}{\gamma+k{u_1^*}}, \quad x^*_3=\frac{\beta}{\gamma}, \quad x_2^*=\frac{g(x_3^*-x_1^*)}{\gamma+k{u_2^*}},
\end{equation}
and the Jacobian matrix:
\begin{equation}
\left.\frac{\partial f(x,u)}{\partial x}\right|_{(x,u)=(x^*,u^*)}=\left[ \begin{array}{ccc}-\gamma-k{u_1^*} & 0 & 0  \\ -g'(x_3^*-x_1^*) & -\gamma-k{u_2^*} & g'(x_3^*-x_1^*) \\ 0  & 0 & -\gamma\end{array}\right]
\end{equation}
has the negative eigenvalues $-\gamma-k{u_1^*}$, $-\gamma-k{u_2^*}$, $-\gamma$, and is thus Hurwitz.  It is clear from (\ref{l1}) and (\ref{l3}) that $x_1(t)$ and $x_3(t)$ converge to $x_1^*$ and $x_3^*$. This means that the first term in (\ref{l2}) converges to $g(x_3^*-x_1^*)$, from which we conclude that $x_2(t)$ converges to $x_2^*$. Thus, $x^*$ is globally asymptotically stable and all other statements in Assumption \ref{characteristic} hold.

To verify Assumption \ref{iomon}, we note that:
\begin{equation}
\frac{\partial f_1}{\partial u_1}=-k{x_1}\le 0, \ \frac{\partial f_2}{\partial u_2}=-k{x_2}\le 0,\ \frac{\partial f_2}{\partial x_1}=-g'(x_3-x_1)\ge 0,\ \frac{\partial f_2}{\partial x_3}=g'(x_3-x_1)\le 0,\ \frac{\partial h_1}{\partial x_2}=1, \ \frac{\partial h_2}{\partial x_1}=1.
\end{equation}
Thus, Lemma \ref{monotonetest} holds with $\delta_1=\delta_2=0$, $\mu_1=\mu_2=1$, $\epsilon_1=\epsilon_2=1$, $\epsilon_3=0$ and, thus, we conclude monotonicity with respect to the orthants $K^U=\mathbb{R}_{\ge 0}^2$, $K^Y=-K^U$, $K^X=\{x\in \mathbb{R}^3 \ | x_1\le 0, x_2\le 0, x_3\ge 0\}$.

\begin{figure}[h]
\vspace{.5cm}
\begin{center}
\mbox{}\setlength{\unitlength}{.5mm}
\begin{picture}(60,70)
\put(-20,20){\psfig{figure=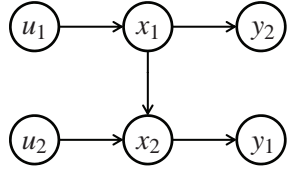,width=100\unitlength}}
 \put(-4,62){$u_1$}\put(-4,31.5){$u_2$}
 \put(26,62){$x_1$}\put(26,31.5){$x_2$}
 \put(57,62){$y_2$}\put(57,31.5){$y_1$}
\end{picture}
\vspace{-1.4cm}
 \caption{\small The incidence graph for system (\ref{l1})-(\ref{l5}), constructed as in Lemma \ref{copy}.}  \label{incidence}
\end{center}
\vspace{-.4cm}
\end{figure}

To show that Assumption \ref{either} holds, we apply the test in Lemma \ref{copy}. The incidence graph for the system (\ref{l1})-(\ref{l5}) restricted to the subset of $\mathscr{X}$ where $x_1>0$, $x_2 > 0$, $x_3=x_3^*$ is as in Figure \ref{incidence}.  From Lemma \ref{copy}, we conclude that the system is excitable by $u_1$, since a directed path connects $u_1$ to both $x_1$ and $x_2$, and transparent from $y_1$, since a directed path connects both $x_1$ and $x_2$ to $y_1$.
\hfill $\Box$

Note that the restriction $x_3=x_3^*$ allowed us to drop $x_3$, which is not excitable by either input, from the incidence graph in Figure \ref{incidence}. Likewise, the restriction $x_j>0$, $j=1,2,$ is critical for the sign-definiteness of $\partial f_j/\partial u_j=-kx_j$, which made it possible to direct an edge from $u_j$ to $x_j$. Because the subset of $\mathscr{X}$ defined by these restrictions is forward invariant and attractive, the $\omega$-limit sets of all solutions starting in $\mathscr{X}$ lie in this subset. Thus, strong monotonicity on this subset, established by Theorem \ref{smon} when the contact graph is bipartite, allows us to conclude generic convergence on $\mathscr{X}$.

We emphasize that the assumption of identical degradation rates for $N$ and $S$ in (\ref{receptor})-(\ref{cofactor}) is essential for the change of coordinates that lead to (\ref{l1})-(\ref{l5}) and that allowed us to conclude monotonicity using Lemma \ref{monotonetest} for orthant orders. It would be interesting to investigate whether monotonicity can be established for nonidentical degradation rates.


\begin{thebibliography}{10}
\providecommand{\url}[1]{#1}
\csname url@samestyle\endcsname
\providecommand{\newblock}{\relax}
\providecommand{\bibinfo}[2]{#2}
\providecommand{\BIBentrySTDinterwordspacing}{\spaceskip=0pt\relax}
\providecommand{\BIBentryALTinterwordstretchfactor}{4}
\providecommand{\BIBentryALTinterwordspacing}{\spaceskip=\fontdimen2\font plus
\BIBentryALTinterwordstretchfactor\fontdimen3\font minus
  \fontdimen4\font\relax}
\providecommand{\BIBforeignlanguage}[2]{{%
\expandafter\ifx\csname l@#1\endcsname\relax
\typeout{** WARNING: IEEEtran.bst: No hyphenation pattern has been}%
\typeout{** loaded for the language `#1'. Using the pattern for}%
\typeout{** the default language instead.}%
\else
\language=\csname l@#1\endcsname
\fi
#2}}
\providecommand{\BIBdecl}{\relax}
\BIBdecl

\bibitem{murray2}
J.~Murray, \emph{Mathematical Biology, II. Spatial Models and Biomedical
  Applications}, 3rd~ed.\hskip 1em plus 0.5em minus 0.4em\relax New York:
  Springer, 2001.

\bibitem{gilbert}
S.~Gilbert, \emph{Developmental Biology}, 9th~ed.\hskip 1em plus 0.5em minus
  0.4em\relax Sinauer Associates, Inc., 2010.

\bibitem{ColMonMaiLew96}
J.~Collier, N.~Monk, P.~Maini, and J.~Lewis, ``Pattern formation by lateral
  inhibition with feedback: a mathematical model of delta-notch intercellular
  signalling,'' \emph{Journal of Theoretical Biology}, vol. 183, pp. 429--446,
  1996.

\bibitem{SprLakLeB10}
D.~Sprinzak, A.~Lakhanpal, L.~LeBon, L.~Santat, M.~Fontes, G.~Anderson,
  J.~Garcia-Ojalvo, and M.~Elowitz, ``Cis-interactions between {Notch} and
  {Delta} generate mutually exclusive signalling states,'' \emph{Nature}, vol.
  465, pp. 86--90, 2010.

\bibitem{SprLakLeB11}
D.~Sprinzak, A.~Lakhanpal, L.~LeBon, J.~Garcia-Ojalvo, and M.~Elowitz, ``Mutual
  inactivation of notch receptors and ligands facilitates developmental
  patterning,'' \emph{{PLoS} Computational Biology}, vol.~7, no.~6, 2011.

\bibitem{AokLow10}
K.~Aoki, E.~Diner, C.~de~Roodenbeke, B.~Burgess, S.~Poole, B.~Braaten,
  A.~Jones, J.~Webb, C.~Hayes, P.~Cotter, and D.~Low, ``A widespread family of
  polymorphic contact-dependent toxin delivery systems in bacteria,''
  \emph{Nature}, vol. 468, pp. 439--442, 2010.

\bibitem{AngSon03}
D.~Angeli and E.~Sontag, ``Monotone control systems,'' \emph{{IEEE}
  Transactions on Automatic Control}, vol.~48, pp. 1684--1698, 2003.

\bibitem{Hirsch-Smith}
M.~Hirsch and H.~Smith, ``Monotone dynamical systems,'' in \emph{Handbook of
  Differential Equations, Ordinary Differential Equations (second
  volume)}.\hskip 1em plus 0.5em minus 0.4em\relax Amsterdam: Elsevier, 2005.

\bibitem{smith}
H.~Smith, \emph{Monotone Dynamical Systems: An Introduction to the Theory of
  Competitive and Cooperative Systems}.\hskip 1em plus 0.5em minus 0.4em\relax
  Providence, RI: American Mathematical Society, 1995.

\bibitem{asratian}
A.~Asratian, T.~Denley, and R.~H{\"{a}}ggkvist, \emph{Bipartite Graphs and
  Their Applications}.\hskip 1em plus 0.5em minus 0.4em\relax Cambridge, UK:
  Cambridge University Press, 1998.

\bibitem{levin}
D.~Levin, Y.~Peres, and E.~Wilmer, \emph{Markov Chains and Mixing Times}.\hskip
  1em plus 0.5em minus 0.4em\relax Providence, Rhode Island: American
  Mathematical Society, 2009.

\bibitem{AngSon04}
D.~Angeli and E.~Sontag, ``Multistability in monotone {I/O} systems,''
  \emph{Systems and Control Letters}, vol.~51, pp. 185--202, 2004.

\bibitem{EncSon05}
G.~Enciso and E.~Sontag, ``Monotone systems under positive feedback:
  multistability and a reduction theorem,'' \emph{Systems and Control Letters},
  vol.~51, pp. 185--202, 2005.

\bibitem{AngSon04b}
D.~Angeli and E.~Sontag, ``Interconnections of monotone systems with
  steady-state characteristics,'' in \emph{Optimal control, stabilization and
  nonsmooth analysis}, M.~de~Queiroz, M.~Malisoff, and P.~Walenski, Eds.\hskip
  1em plus 0.5em minus 0.4em\relax Springer, 2004, vol. 301, pp. 135--154.

\bibitem{berman}
A.~Berman and R.~Plemmons, \emph{Nonnegative Matrices in the Mathematical
  Sciences}.\hskip 1em plus 0.5em minus 0.4em\relax Philadelphia: Society for
  Industrial and Applied Mathematics, Classics in Applied Mathematics, 1994,
  originally published by Academic Press, New York, 1979.

\bibitem{kocic}
V.~Kocic and G.~Ladas, \emph{Global Behavior of Nonlinear Difference Equations
  of Higher Order with Applications}.\hskip 1em plus 0.5em minus 0.4em\relax
  Kluwer Academic Publishers, 1993.

\bibitem{WeaShe01}
H.~Wearing and J.~Sherrat, ``Analysis of juxtacrine patterns,'' \emph{SIAM
  Journal on Applied Mathematics}, vol.~62, pp. 283--309, 2001.

\bibitem{Pla01}
E.~Plahte, ``Pattern formation in discrete cell lattices,'' \emph{Journal of
  Mathematical Biology}, vol.~43, pp. 411--445, 2001.

\end{thebibliography}
\end{document}